# A unifying framework for *k*-statistics, polykays and their multivariate generalizations

ELVIRA DI NARDO[1,*], GIUSEPPE GUARINO[2] and DOMENICO SENATO[1,**]

[1]*Dipartimento di Matematica e Informatica, Università degli Studi della Basilicata, Viale dell'Ateneo Lucano 10, I-85100 Potenza, Italy.* E-mail: [*]*elvira.dinardo@unibas.it*; [**]*domenico.senato@unibas.it*

[2]*Medical School, Università Cattolica del Sacro Cuore (Rome branch), Largo Agostino Gemelli 8, I-00168 Roma, Italy.* E-mail: *giuseppe.guarino@rete.basilicata.it*

Through the classical umbral calculus, we provide a unifying syntax for single and multivariate *k*-statistics, polykays and multivariate polykays. From a combinatorial point of view, we revisit the theory as exposed by Stuart and Ord, taking into account the Doubilet approach to symmetric functions. Moreover, by using exponential polynomials rather than set partitions, we provide a new formula for *k*-statistics that results in a very fast algorithm to generate such estimators.

## 1. Introduction

The theory of *k*-statistics has a long history. Here, we outline the literature, beginning with Fisher (1929), who rediscovered the half-invariants theory of Thiele (1897). He introduced *k*-statistics (single and multivariate) as new symmetric functions of a random sample, aiming to estimate cumulants without using moment estimators. Dressel (1940) developed a theory of more general functions, later resumed by Tukey (1950), who named them *polykays*. Both Tukey (1956) and Wishart (1952) developed methods to express polykays in terms of Fisher's *k*-statistics. These methods are straightforward enough, but their execution leads to intricate computations and some cumbersome expressions, except in very simple cases. Later, many authors tried to simplify the matter. Kaplan (1952) resorted to tensor notation in order to simplify multivariate *k*-statistics. Good (1975) gave an interpretation of cumulants as coefficients of the Fourier transform of the randomly ordered sample and used this formula in order to obtain expressions for single *k*-statistics (1977). The whole subject was later described in great detail by Stuart and Ord (1987).

In the 1980's, tensor notation was employed by Speed (1983, 1986a, 1986b, 1986c), Speed and Silcock (1988a, 1988b) and extended to polykays and single $k$-statistics. This extension reveals the coefficients defining polykays to be values of the Möbius function over the lattice of set partitions. As a consequence, Speed used the set theoretic approach to symmetric functions introduced by Doubilet (1972). In the same period, McCullagh (1984, 1987) simplified the tensor notation of Kaplan by introducing the notion of generalized cumulants. Symbolic operators for expectation and the derivation of unbiased estimates for multiple sums were introduced by Andrews (2001), Andrews and Stafford (1998, 2000). Algorithms to compute $k$-statistics and their generalizations were derived from such techniques.

In 1994, Rota and Taylor brought new life to the celebrated umbral calculus, with the intention of restoring (within a modern formal setting) the computational power originally dealt with in the writings of Blissard, Cayley and Sylvester. Their basic device was the representation of a unital sequence of numbers by a symbol, called the *umbra*, via an operator resembling the expectation operator of random variables. This shifts the attention to the moment generating function rather than the probability generating one, with the latter usually linking the umbral calculus of Roman and Rota (1978) to probability theory (for a survey, see Di Bucchianico (1997)). The umbral calculus of Rota and Taylor, known as the *classical umbral calculus*, has been developed by Di Nardo and Senato (2001, 2006a), paying particular attention to relations with probability theory. It has also been successfully applied in wavelet theory (Saliani and Senato (2006) and Shen (1999)).

In this paper, we show how the classical umbral calculus provides a unifying framework for $k$-statistics and their generalizations. Most of the results can be found in some form in the literature; nevertheless, we feel that the umbral approach unifies, generalizes and simplifies their presentation. In order to demonstrate the power of the umbral methodologies, we conclude the paper by presenting an algorithm that generates $k$-statistics in a very short computational time compared with existing procedures.

Section 2 is aimed at readers unaware of the classical umbral calculus. Here, we recall basic definitions and terminology. In Section 3, we consider symmetric polynomials in umbral frameworks and introduce some combinatorial tools. If the indeterminates of symmetric polynomials are replaced by uncorrelated umbrae, surprisingly compact expressions for both symmetric polynomials and relations between bases can be achieved. From a statistical point of view, this means obtaining elementary formulae connecting symmetric polynomials in the data points. Thus, the umbral approach recovers the methods exposed by Stuart and Ord, whose main limitation is the complexity of both the expressions and the procedures. In Section 4, the theory of $k$-statistics and polykays is completely rewritten, stressing the power of the umbral methods by means of some examples. In Section 5, we introduce the notion of umbra indexed by a multiset, in order to get umbral expressions for multivariate $k$-statistics and multivariate polykays. Through this device, we represent multivariate moments and multivariate cumulants by umbrae. Moreover, we state identities involving multivariate moments and multivariate cumulants by simply characterizing a suitable multiset indexing. In the last section, we give a very fast algorithm to compute $k$-statistics, originating from a new umbral expression for such estimators. Comparisons of computational times achieved by the Andrews and Stafford algorithm, that of MATHSTATICA (Rose and Smith (2002)) and the proposed umbral algorithm end this paper.

## 2. Background to classical umbral calculus

This section reviews notation and terminology useful when handling umbrae. More details and technicalities can be found in Di Nardo and Senato (2001, 2006a).

Formally, umbral calculus is a syntax consisting of the following data:

(i) a set $A = \{\alpha, \beta, \ldots\}$, called the *alphabet*, whose elements are named *umbrae*;

(ii) a commutative integral domain $R$ whose quotient field is of characteristic zero;

(iii) a linear functional $E$, called the *evaluation*, defined on the polynomial ring $R[A]$ and taking values in $R$ such that

- $E[1] = 1$;
- $E[\alpha^i \beta^j \cdots \gamma^k] = E[\alpha^i] E[\beta^j] \cdots E[\gamma^k]$ for any set of distinct umbrae in $A$ and for $i, j, \ldots, k$ non-negative integers (*uncorrelation property*);

(iv) an element $\varepsilon \in A$, called an *augmentation*, such that $E[\varepsilon^n] = 0$ for every $n \geq 1$;

(v) an element $u \in A$, called a *unity* umbra, such that $E[u^n] = 1$ for every $n \geq 1$.

An *umbral polynomial* is a polynomial $p \in R[A]$. The support of $p$ is the set of all umbrae occurring in $p$. If $p$ and $q$ are two umbral polynomials, then

(i) $p$ and $q$ are *uncorrelated* if and only if their supports are disjoint;

(ii) $p$ and $q$ are *umbrally equivalent* iff $E[p] = E[q]$ (in symbols, $p \simeq q$).

The *moments* of an umbra $\alpha$ are the elements $a_n \in R$ such that $E[\alpha^n] = a_n$ for $n \geq 0$ and we say that the umbra $\alpha$ *represents* the sequence of moments $1, a_1, a_2, \ldots$.

It is possible that two distinct umbrae represent the same sequence of moments. In such a case, they are called *similar umbrae*. More formally, two umbrae $\alpha$ and $\gamma$ are said to be similar when

$$E[\alpha^n] = E[\gamma^n] \qquad \text{for every } n \geq 0 \qquad (\text{in symbols}, \alpha \equiv \gamma).$$

Furthermore, given a sequence $1, a_1, a_2, \ldots$ in $R$, there are infinitely many distinct, and thus similar, umbrae representing the sequence.

The *factorial moments* of an umbra $\alpha$ are the elements $a_{(n)} \in R$ corresponding to umbral polynomials $(\alpha)_n = \alpha(\alpha - 1) \cdots (\alpha - n + 1), n \geq 1$, via the evaluation $E$, that is, $E[(\alpha)_n] = a_{(n)}$.

***Example 2.1 (Singleton umbra).*** The *singleton umbra* $\chi$ is the umbra whose moments are all zero, except the first $E[\chi] = 1$. Its factorial moments are $x_{(n)} = (-1)^{n-1}(n-1)!$.

***Example 2.2 (Bell umbra).*** The *Bell umbra* $\beta$ is the umbra whose factorial moments are all equal to 1, that is, $E[(\beta)_n] = 1$ for every $n \geq 1$. Its moments are the Bell numbers.

Thanks to the notion of similar umbrae, it is possible to extend the alphabet $A$ with the so-called *auxiliary umbrae* obtained via operations among similar umbrae. This leads to the construction of a *saturated umbral calculus* in which auxiliary umbrae are treated as elements of the alphabet (Rota and Taylor (1994)). Auxiliary umbrae come into play when dealing with products of moments. Note that $a_i a_j \neq E[\alpha^i \alpha^j]$, with $a_i = E[\alpha^i]$ and $a_j = E[\alpha^j]$, but $a_i a_j = E[\alpha^i \alpha'^j]$, with $\alpha \equiv \alpha'$ and $\alpha'$ uncorrelated with $\alpha$. As a consequence, the umbra $\alpha + \alpha'$ represents the sequence

$$\sum_{k=0}^{n} \binom{n}{k} a_{n-k} a_k.$$

We may then denote the umbra similar to $\alpha + \alpha'$ by the auxiliary umbra $2.\alpha$. In a saturated umbral calculus, the umbra $2.\alpha$ is considered as an element of the alphabet $A$.

In the following, we focus attention on auxiliary umbrae which will subsequently have a special role. Let $\{\alpha_1, \alpha_2, \ldots, \alpha_n\}$ be a set of $n$ uncorrelated umbrae similar to an umbra $\alpha$. The symbol $n.\alpha$ denotes an auxiliary umbra similar to the sum $\alpha_1 + \alpha_2 + \cdots + \alpha_n$. The symbol $\alpha^{\cdot n}$ denotes an auxiliary umbra similar to the product $\alpha_1 \alpha_2 \cdots \alpha_n$. Properties of these auxiliary umbrae are extensively described in Di Nardo and Senato (2001) and will be recalled whenever necessary.

*Remark 2.1.* If $n \neq m$, then $n.\alpha$ is uncorrelated with $m.\alpha$ and $\alpha^{\cdot n}$ is uncorrelated with $\alpha^{\cdot m}$. If $p$ and $q$ are correlated umbral polynomials, then $n.p \simeq p_1 + \cdots + p_n$ is correlated with $n.q \simeq q_1 + \cdots + q_n$, and $p_i$ is correlated with $q_i$, but uncorrelated with $q_j$ if $i \neq j$.

If the umbra $\alpha$ represents the sequence $1, a_1, a_2, \ldots$, then $E[(\alpha^{\cdot n})^k] = a_k^n$ for non-negative integers $k$ and $n$.

In Di Nardo and Senato (2006a), the moments of $n.\alpha$ have been expressed in terms of moments of $\alpha$, by means of the Bell exponential polynomials. Here, we adopt a different point of view referring to the notion of integer partitions, obviously connected to the Bell exponential polynomials through well-known relations.

Recall that a partition of an integer $i$ is a sequence $\lambda = (\lambda_1, \lambda_2, \ldots, \lambda_t)$, where $\lambda_j$ are weakly decreasing integers and $\sum_{j=1}^{t} \lambda_j = i$. The integers $\lambda_j$ are named *parts* of $\lambda$. The *length* of $\lambda$ is the number of its parts and will be indicated by $\nu_\lambda$. A different notation is $\lambda = (1^{r_1}, 2^{r_2}, \ldots)$, where $r_j$ is the number of parts of $\lambda$ equal to $j$ and $r_1 + r_2 + \cdots = \nu_\lambda$. We use the classical notation $\lambda \vdash i$ to denote that $\lambda$ is a partition of $i$. Through the multinomial expansion theorem, powers of $n.\alpha$ can be seen to be umbrally equivalent to the umbral polynomials

$$(n.\alpha)^i \simeq \sum_{\lambda \vdash i} (n)_{\nu_\lambda} d_\lambda \alpha_\lambda, \tag{2.1}$$

where the sum is over all partitions of the integer $i$, $(n)_{\nu_\lambda} = 0$ when $\nu_\lambda > n$,

$$d_\lambda = \frac{i!}{r_1! r_2! \cdots} \frac{1}{(1!)^{r_1} (2!)^{r_2} \cdots} \quad \text{and} \quad \alpha_\lambda \equiv (\alpha_{j_1})^{\cdot r_1} (\alpha_{j_2}^2)^{\cdot r_2} \cdots, \tag{2.2}$$

with $\{j_i\}$ distinct integers chosen in $\{1, 2, \ldots, n\} = [n]$. In particular, when evaluating the powers of $n.\alpha$ in (2.1), we have

$$E[(n.\alpha)^i] = \sum_{\lambda \vdash i} (n)_{\nu_\lambda} d_\lambda a_\lambda, \tag{2.3}$$

where $a_\lambda = a_1^{r_1} a_2^{r_2} \cdots$.

A feature of classical umbral calculus is the construction of new auxiliary umbrae by symbolic substitution. For example, if we replace the integer $n$ in $n.\alpha$ with an umbra $\gamma$, then, from (2.1), the new auxiliary umbra $\gamma.\alpha$ has powers

$$(\gamma.\alpha)^i \simeq \sum_{\lambda \vdash i} (\gamma)_{\nu_\lambda} d_\lambda \alpha_\lambda. \tag{2.4}$$

Equivalence (2.4) has been formally proven by using the notion of the generating function of an umbra; for further details, see Di Nardo and Senato (2001).

In the dot product $\gamma.\alpha$, replacing the umbra $\gamma$ with the umbra $\gamma.\beta$, we obtain the composition umbra of $\alpha$ and $\gamma$, that is, $\gamma.\beta.\alpha$. Its powers are

$$(\gamma.\beta.\alpha)^i \simeq \sum_{\lambda \vdash i} \gamma^{\nu_\lambda} d_\lambda \alpha_\lambda. \tag{2.5}$$

The *compositional inverse* of an umbra $\alpha$ is the umbra $\alpha^{<-1>}$ such that

$$\alpha^{<-1>}.\beta.\alpha \equiv \alpha.\beta.\alpha^{<-1>} \equiv \chi.$$

The compositional inverse of an umbra was introduced in Di Nardo and Senato (2001) in order to invert exponential power series. It has also been used to give a simple proof of the Lagrange inversion formula in umbral terms. Moreover, the compositional inverse of an umbra provides a link between the Bell and singleton umbra. For the purpose of this paper, we consider the compositional inverse of the unity umbra.

In the following examples, powers of fundamental auxiliary umbrae are given via (2.4). Properties of such umbrae are described in Di Nardo and Senato (2001, 2006a).

*Example 2.3 ($\alpha$-partition umbra).* The umbra $\beta.\alpha$, with $\beta$ the Bell umbra, is called the *$\alpha$-partition umbra*. By virtue of (2.4), its powers are

$$(\beta.\alpha)^i \simeq \sum_{\lambda \vdash i} d_\lambda \alpha_\lambda \tag{2.6}$$

because the factorial moments of the Bell umbra are all equal to 1 (see Example 2.2). In particular, we have

$$\beta.u^{<-1>} \equiv \chi, \qquad \beta.u \equiv \beta, \qquad \beta.\chi \equiv u, \tag{2.7}$$

where $u^{<-1>}$ denotes the compositional inverse of $u$.

***Example 2.4 (α-cumulant umbra).*** The umbra $\chi \boldsymbol{.} \alpha$, with $\chi$ the singleton umbra, is called the *α-cumulant umbra*. By virtue of (2.4), its powers are

$$(\chi \boldsymbol{.} \alpha)^i \simeq \sum_{\lambda \vdash i} x_{(\nu_\lambda)} d_\lambda \alpha_\lambda, \tag{2.8}$$

where $x_{(\nu_\lambda)}$ are the factorial moments of the umbra $\chi$ (see Example 2.1). In particular, it is possible to prove that

$$\chi \boldsymbol{.} \beta \equiv u, \qquad \chi \boldsymbol{.} \chi \equiv u^{<-1>}.$$

***Example 2.5 (α-factorial umbra).*** The umbra $\alpha \boldsymbol{.} \chi$ is called the *α-factorial umbra* and its moments are the factorial moments of $\alpha$, that is, $(\alpha \boldsymbol{.} \chi)^i \simeq (\alpha)_i$.

The disjoint sum of $\alpha$ and $\gamma$ is the umbra whose moments are the sum of $n$th moments of $\alpha$ and $\gamma$ respectively (Di Nardo and Senato (2006a))—in symbols,

$$(\alpha \dotplus \gamma)^n \simeq \alpha^n + \gamma^n \qquad \text{for every } n > 0.$$

For instance, it is possible to prove

$$\chi \boldsymbol{.} \alpha \dotplus \chi \boldsymbol{.} \gamma \equiv \chi \boldsymbol{.} (\alpha + \gamma), \tag{2.9}$$

the well-known additive property of cumulants. In the following, we denote by $\dotplus_{i=1}^{n} \alpha_i$ the disjoint sum of $n$ umbrae and by $\dotplus_n \alpha$ the disjoint sum of $n$ times the umbra $\alpha$.

## 3. Symmetric polynomial umbrae

We begin by recalling the definitions of the four classical bases of the algebra of symmetric polynomials in the variables $x_1, x_2, \ldots, x_n$. They are:
*elementary symmetric polynomials*

$$e_k = \sum_{1 \le j_1 < j_2 < \cdots < j_k \le n} x_{j_1} x_{j_2} \cdots x_{j_k};$$

*power sum symmetric polynomials*

$$s_r = \sum_{i=1}^{n} x_i^r, \qquad r = 1, 2, \ldots;$$

*monomial symmetric polynomials*

$$m_\lambda = \sum x_1^{\lambda_1} \cdots x_t^{\lambda_t},$$

where the sum is over all distinct monomials having exponents $\lambda_1, \ldots, \lambda_t$;

*complete homogeneous symmetric polynomials*

$$h_i = \sum_{\lambda \vdash i} m_\lambda = \sum_{1 \leq j_1 \leq j_2 \leq \cdots \leq j_k \leq n} x_{j_1} x_{j_2} \cdots x_{j_k}.$$

In the following, we replace the commutative integral domain $R$ by $K[x_1, x_2, \ldots, x_n]$, where $K$ is a field of characteristic zero and $x_1, x_2, \ldots, x_n$ are variables. Therefore, the uncorrelation property (iii) of Section 2 must be rewritten as

$$E[1] = 1; \qquad E[x_i x_j \cdots \alpha^k \beta^l \cdots] = x_i x_j \cdots E[\alpha^k] E[\beta^l] \cdots$$

for any set of distinct umbrae in $A$, for $i, j, \ldots \in [n]$ and for non-negative integers $k, l, \ldots$.

In $K[x_1, x_2, \ldots, x_n][A]$, an umbra is said to be a *scalar umbra* when its moments are elements of $K$, while it is said to be a *polynomial umbra* if its moments are polynomials of $K[x_1, x_2, \ldots, x_n]$. A polynomial umbra is said to be *symmetric* when its moments are symmetric polynomials in $K[x_1, x_2, \ldots, x_n]$.

A sequence of polynomials $p_0, p_1, \ldots \in K[x_1, x_2, \ldots, x_n]$ is umbrally represented by a polynomial umbra if $p_0 = 1$ and $p_n$ is of degree $n$ for every nonnegative integer $n$. The four classical bases of the algebra of symmetric polynomials are all represented by symmetric polynomial umbrae. In particular, we call the polynomial umbra $\epsilon$ such that

$$E[\epsilon^k] = \begin{cases} e_k, & k = 1, 2, \ldots, n, \\ 0, & k = n+1, n+2, \ldots \end{cases}$$

elementary polynomial umbra. We call the polynomial umbra $\sigma$ such that $E[\sigma^r] = s_r$ for every nonnegative integer $r$, power sum polynomial umbra.

**Proposition 3.1 (Elementary polynomial umbra).** *If $\chi_1, \ldots, \chi_n$ are $n$ uncorrelated umbrae similar to the singleton umbra, then*

$$\epsilon^k \simeq \frac{(\chi_1 x_1 + \cdots + \chi_n x_n)^k}{k!}, \qquad k = 1, 2, \ldots, \tag{3.1}$$

*where $\epsilon$ is the elementary polynomial umbra.*

**Proof.** For $k = 1, \ldots, n$, the result follows by applying the evaluation $E$ to the multinomial expansion of $(\chi_1 x_1 + \cdots + \chi_n x_n)^k$. There are vanishing terms corresponding to the powers of $\chi$ greater than 1. Only $k!$ monomials of the form $\chi_{j_1} x_{j_1} \chi_{j_2} x_{j_2} \cdots \chi_{j_k} x_{j_k}$ have a non-zero evaluation. Instead, for $k = n+1, n+2, \ldots$, the evaluation $E$ gives zero since at least one power of $\chi$ greater than 1 occurs in each monomial of the multinomial expansion. $\square$

**Proposition 3.2 (Power sum polynomial umbra).** *If $u$ is the unity umbra and $\sigma$ is the power sum polynomial umbra, then*

$$\sigma \equiv (\dot{+}_{i=1}^{n} \, u x_i).$$

The following theorem gives an umbral relation between the elementary polynomial umbra and the power sum polynomial umbra.

**Theorem 3.1.** *If $\chi_1, \ldots, \chi_n$ are $n$ uncorrelated umbrae similar to the singleton umbra $\chi$ and $\sigma$ is the power sum polynomial umbra, then*

$$\chi \boldsymbol{\cdot} (\chi_1 x_1 + \cdots + \chi_n x_n) \equiv (\chi \boldsymbol{\cdot} \chi) \sigma. \tag{3.2}$$

**Proof.** Equivalence (3.2) follows from (2.9), observing that

$$\chi \boldsymbol{\cdot} (\chi_1 x_1 + \cdots + \chi_n x_n) \equiv \dot{+}_{i=1}^n \chi \boldsymbol{\cdot} (\chi_i x_i) \equiv \dot{+}_{i=1}^n (\chi \boldsymbol{\cdot} \chi) x_i \equiv (\chi \boldsymbol{\cdot} \chi) \sigma. \tag{3.3}$$

$\square$

The next corollary points out a deeper meaning of Theorem 3.1, that is, by evaluating the moments of the umbrae in (3.2), the well-known relations between power sum symmetric polynomials $\{s_r\}$ and elementary symmetric polynomials $\{e_k\}$ are recovered.

**Corollary 3.1.** *If $\{s_r\}$ are the power symmetric polynomials and $\{e_k\}$ are the elementary symmetric polynomials, then*

$$\frac{(-1)^{i-1}}{i} s_i = \sum_{\lambda \vdash i} x_{(\nu_\lambda)} \frac{e_1^{r_1}}{r_1!} \frac{e_2^{r_2}}{r_2!} \cdots \qquad e_i = \frac{1}{i!} \sum_{\lambda \vdash i} d_\lambda (x_{(1)} s_1)^{r_1} (x_{(2)} s_2)^{r_2} \cdots. \tag{3.4}$$

**Proof.** From (3.1), we have

$$E[(\chi_1 x_1 + \cdots + \chi_n x_n)_\lambda] = E[\chi_1 x_1 + \cdots + \chi_n x_n]^{r_1} E[(\chi_1 x_1 + \cdots + \chi_n x_n)^2]^{r_2} \cdots$$
$$= (1!)^{r_1} e_1^{r_1} (2!)^{r_2} e_2^{r_2} \cdots$$

and so

$$E[(\chi \boldsymbol{\cdot} (\chi_1 x_1 + \cdots + \chi_n x_n))^i] = i! \sum_{\lambda \vdash i} x_{(\nu_\lambda)} \frac{e_1^{r_1}}{r_1!} \frac{e_2^{r_2}}{r_2!} \cdots.$$

By Theorem 3.1, we have

$$E[(\chi \boldsymbol{\cdot} (\chi_1 x_1 + \cdots + \chi_n x_n))^i] = E[(\chi \boldsymbol{\cdot} \chi)^i \sigma^i]$$

and the former identity in (3.4) follows, observing that

$$E[(\chi \boldsymbol{\cdot} \chi)^i \sigma^i] = (-1)^{i-1}(i-1)! s_i.$$

Taking the right dot product with the Bell umbra $\beta$ in (3.2), we have

$$\beta \boldsymbol{\cdot} [(\chi \boldsymbol{\cdot} \chi) \sigma] \equiv \beta \boldsymbol{\cdot} \chi \boldsymbol{\cdot} (\chi_1 x_1 + \cdots + \chi_n x_n)$$

and, by virtue of the third equivalence in (2.7), we have

$$\beta.[(\chi.\chi)\sigma] \equiv \chi_1 x_1 + \cdots + \chi_n x_n.$$

Equivalence (2.6) gives

$$(\beta.[(\chi.\chi)\sigma])^i \simeq \sum_{\lambda \vdash i} d_\lambda (\chi.\chi)_\lambda \sigma_\lambda \simeq i! \epsilon^i.$$

The latter identity in (3.4) follows by observing that $E[\sigma_\lambda] = (s_1)^{r_1}(s_2)^{r_2}\cdots$ and $E[(\chi.\chi)_\lambda] = (x_{(1)})^{r_1}(x_{(2)})^{r_2}\cdots$. $\square$

The umbral expression of $m_\lambda$ requires the introduction of augmented monomial symmetric polynomials $\tilde{m}_\lambda$. Let $\lambda = (1^{r_1}, 2^{r_2}, \ldots)$ be a partition of the integer $i \leq n$. Augmented monomial symmetric polynomials are defined as

$$\tilde{m}_\lambda = \sum_{j_1 \neq \cdots \neq j_{r_1} \neq j_{r_1+1} \neq \cdots \neq j_{r_1+r_2} \neq \cdots} x_{j_1} \cdots x_{j_{r_1}} x_{j_{r_1+1}}^2 \cdots x_{j_{r_1+r_2}}^2 \cdots.$$

The next proposition can be proven using the same approach as used in the proof of Proposition 3.1.

**Proposition 3.3.** *If $\chi_1, \ldots, \chi_n$ are $n$ uncorrelated umbrae similar to the singleton umbra, then*

$$\tilde{m}_\lambda \simeq (\chi_1 x_1 + \cdots + \chi_n x_n)^{r_1} (\chi_1 x_1^2 + \cdots + \chi_n x_n^2)^{r_2} \cdots. \tag{3.5}$$

The next corollary follows by recalling that $m_\lambda = \tilde{m}_\lambda/(r_1! r_2! \cdots)$.

**Corollary 3.2.**

$$m_\lambda \simeq \frac{(\chi_1 x_1 + \cdots + \chi_n x_n)^{r_1}}{r_1!} \frac{(\chi_1 x_1^2 + \cdots + \chi_n x_n^2)^{r_2}}{r_2!} \cdots.$$

In order to characterize the symmetric polynomial umbra representing the complete homogeneous symmetric polynomials, we need to recall the notion of the *inverse* of an umbra. Two umbrae $\alpha$ and $\gamma$ are said to be *inverse* to each other when $\alpha + \gamma \equiv \varepsilon$. The inverse of the umbra $\alpha$ is denoted by $-1.\alpha$. Note that, in dealing with a saturated umbral calculus, the inverse of an umbra is not unique, but any two inverse umbrae of the same umbra are similar.

**Proposition 3.4 (Complete homogeneous polynomial umbra).** *If $\chi_1, \ldots, \chi_n$ are $n$ uncorrelated umbrae similar to the singleton umbra, then*

$$h_i \simeq \frac{\{-1.[\chi_1(-x_1) + \cdots + \chi_n(-x_n)]\}^i}{i!}, \qquad i = 1, 2, \ldots. \tag{3.6}$$

**Proof.** From the multinomial expansion theorem, we have

$$[-1 \boldsymbol{.} (-\chi_1 x_1) + \cdots + -1 \boldsymbol{.} (-\chi_n x_n)]^i$$
$$\simeq i! \sum_{|\lambda|=i} [-1 \boldsymbol{.} (-\chi_{j_1})]^{\boldsymbol{.} r_1} ([-1 \boldsymbol{.} (-\chi_{j_2})]^2)^{\boldsymbol{.} r_2} \cdots \frac{\tilde{m}_\lambda}{(1!)^{r_1} r_1! (2!)^{r_2} r_2! \cdots},$$

where $j_1, j_2, \ldots$ are distinct integers chosen in $[n]$. As $([-1 \boldsymbol{.} (-\chi)]^i)^{\boldsymbol{.} r_i} \simeq (i!)^{r_i}$, the result follows from Corollary 3.2 and Proposition 3.3. □

Note that equivalences (3.1) and (3.6) are umbral versions of the well-known identities

$$\sum_k e_k t^k = \prod_{i=1}^n (1 + x_i t), \qquad \sum_k h_k t^k = \frac{1}{\prod_{i=1}^n (1 - x_i t)}.$$

**Proposition 3.5.** *If $\chi_1, \ldots, \chi_n$ are $n$ uncorrelated umbrae similar to the singleton umbra $\chi$ and $\sigma$ is the power sum polynomial umbra, then*

$$-\chi \boldsymbol{.} [\chi_1(-x_1) + \cdots + \chi_n(-x_n)] \equiv [-\chi \boldsymbol{.} (-\chi)] \sigma. \tag{3.7}$$

**Proof.** Equivalence (3.7) follows by replacing the umbra $\chi$ with $-\chi$ and the umbra $u^{<-1>} \equiv \chi \boldsymbol{.} \chi$ with $(-u)^{<-1>} \equiv (-\chi) \boldsymbol{.} (-\chi)$ in (3.3). □

Equivalence (3.7) is an umbral version of the well-known relations between power sum symmetric polynomials $\{s_r\}$ and complete homogeneous symmetric polynomials $\{h_k\}$.

### 3.1. Umbral symmetric polynomials

Assume that we replace the variables $x_1, x_2, \ldots, x_n$ in the umbral polynomials $\chi_1 x_1 + \cdots + \chi_n x_n$ with $n$ uncorrelated umbrae $\alpha_1, \alpha_2, \ldots, \alpha_n$ similar to an umbra $\alpha$. Since

$$\chi_1 \alpha_1 + \cdots + \chi_n \alpha_n \equiv n \boldsymbol{.} (\chi \alpha),$$

from Proposition 3.1, we have

$$\frac{[n \boldsymbol{.} (\chi \alpha)]^k}{k!} \simeq e_k(\alpha_1, \ldots, \alpha_n),$$

where $e_k(\alpha_1, \ldots, \alpha_n)$ are umbral elementary symmetric polynomials in $K[A]$. The same substitution in $u x_1 \dot{+} \cdots \dot{+} u x_n$ gives

$$\dot{+}_{i=1}^n u \alpha_i \equiv \dot{+}_n \alpha \Rightarrow (\dot{+}_n \alpha)^k \simeq \alpha_1^k + \cdots + \alpha_n^k \equiv n \boldsymbol{.} \alpha^k$$

so that $n \boldsymbol{.} \alpha^k \equiv s_k(\alpha_1, \ldots, \alpha_n)$, where $s_k(\alpha_1, \ldots, \alpha_n)$ are umbral power sum symmetric polynomials in $K[A]$. By using the latter identity in (3.4), we prove the following proposition.

**Proposition 3.6.** *If $\chi$ is the singleton umbra and $\alpha \in A$, then*

$$[n.(\chi\alpha)]^k \simeq \sum_{\lambda \vdash k} d_\lambda (\chi.\chi)_\lambda (n.\alpha)^{r_1} (n.\alpha^2)^{r_2} \cdots. \tag{3.8}$$

Due to (3.5), umbral augmented symmetric polynomials $\tilde{m}_\lambda$ take the form

$$\tilde{m}_\lambda \simeq [n.(\chi\alpha)]^{r_1} [n.(\chi\alpha^2)]^{r_2} \cdots.$$

**Theorem 3.2.** *If $\lambda \vdash i$, then*

$$[n.(\chi\alpha)]^{r_1} [n.(\chi\alpha^2)]^{r_2} \cdots \simeq (n)_{\nu_\lambda} \alpha_\lambda. \tag{3.9}$$

(For the proof, see Di Nardo and Senato (2006b).)

Statistically speaking, the last theorem states how to estimate products of moments by using only $n$ sampled items. Moreover, by using Theorem 3.2 and (2.1), we are able to give a sort of inversion formula of equivalence (3.8), that is,

$$(n.\alpha)^k \simeq \sum_{\lambda \vdash k} d_\lambda [n.(\chi\alpha)]^{r_1} [n.(\chi\alpha^2)]^{r_2} \cdots. \tag{3.10}$$

Equivalences (3.8) and (3.10) can be rewritten by using set partitions instead of integer partitions. Afterward, the use of set partitions allows a natural generalization of these equivalences.

Let $C$ be a subset of $K[A]$ such that $|C| = n$. Recall that a partition $\pi$ of $C$ is a collection $\pi = \{B_1, B_2, \ldots, B_k\}$ with $k \leq n$ disjoint and non-empty subsets of $C$ whose union is $C$. We denote by $\Pi_n$ the set of all partitions of $C$.

Let $\{\alpha_1, \alpha_2, \ldots, \alpha_n\}$ be a set of $n$ uncorrelated umbrae similar to an umbra $\alpha$. We will denote by $\alpha^{.\pi}$ the umbra

$$\alpha^{.\pi} \equiv \alpha_{i_1}^{|B_1|} \alpha_{i_2}^{|B_2|} \cdots \alpha_{i_k}^{|B_k|}, \tag{3.11}$$

where $\pi = \{B_1, B_2, \ldots, B_k\}$ is a partition of $\{\alpha_1, \alpha_2, \ldots, \alpha_n\}$ and $i_1, i_2, \ldots, i_k$ are distinct integers chosen in $[n]$. In particular, $\alpha^{.\pi} \equiv \alpha_\lambda$, where $\lambda$ is the partition of the integer $n$ determined by $\pi$. Indeed, a set partition is said to be of type $\lambda = (1^{r_1}, 2^{r_2}, \ldots)$ if there are $r_1$ blocks of cardinality 1, $r_2$ blocks of cardinality 2 and so on. The number of set partitions of type $\lambda$ is $d_\lambda$, as given in (2.2).

**Proposition 3.7.** *If $\Pi_k$ is the set of all partitions of $[k]$ and $\alpha \in A$, then*

$$(n.\alpha)^k \simeq \sum_{\pi \in \Pi_k} [n.(\chi\alpha)]^{r_1} [n.(\chi\alpha^2)]^{r_2} \cdots, \tag{3.12}$$

$$[n.(\chi\alpha)]^k \simeq \sum_{\pi \in \Pi_k} (\chi.\chi)^{.\pi} (n.\alpha)^{r_1} (n.\alpha^2)^{r_2} \cdots. \tag{3.13}$$

**Proof.** Equivalence (3.12) follows directly from (3.10). From (3.11) we have $(\chi \cdot \chi)_\lambda \equiv (\chi \cdot \chi)^{\cdot \pi}$. Hence equivalence (3.8) imples (3.13). □

In order to write products like $(n.\alpha)^{r_1}(n.\alpha^2)^{r_2} \cdots$ using a single symbol such as $\alpha_\lambda$ or $\alpha^{\cdot \pi}$, we need the notion of the multiset. This will be introduced in the next section.

### 3.2. Some necessary combinatorics

A *multiset* $M$ is a pair $(\bar{M}, f)$, where $\bar{M}$ is a set, called the *support* of the multiset, and $f$ is a function from $\bar{M}$ to the non-negative integers. For each $\mu \in \bar{M}$, $f(\mu)$ is called the *multiplicity* of $\mu$. The *length* of the multiset $(\bar{M}, f)$, usually denoted by $|M|$, is the sum of multiplicities of all elements of $\bar{M}$, that is,

$$|M| = \sum_{\mu \in \bar{M}} f(\mu).$$

From now on, we denote a multiset $(\bar{M}, f)$ simply by $M$. A multiset $M_i = (\bar{M}_i, f_i)$ is called a *submultiset* of $M = (\bar{M}, f)$ if $\bar{M}_i \subseteq \bar{M}$ and $f_i(\mu) \leq f(\mu)$ for every $\mu \in \bar{M}_i$.

Let $M$ be a multiset of umbral monomials. When the support of $M$ is a finite set, say $\bar{M} = \{\mu_1, \mu_2, \ldots, \mu_k\}$, we will write

$$M = \{\mu_1^{(f(\mu_1))}, \mu_2^{(f(\mu_2))}, \ldots, \mu_k^{(f(\mu_k))}\} \quad \text{or} \quad M = \{\underbrace{\mu_1, \ldots, \mu_1}_{f(\mu_1)}, \ldots, \underbrace{\mu_k, \ldots, \mu_k}_{f(\mu_k)}\}.$$

Set

$$\mu_M = \prod_{\mu \in \bar{M}} \mu^{f(\mu)}. \tag{3.14}$$

When $f(\mu) = 1$ for every $\mu \in \bar{M}$, the multiset $M$ is simply a set $B$ of umbral monomials in $K[A]$ and hence (3.14) becomes

$$\mu_B = \prod_{\mu \in B} \mu.$$

If $B = \{\alpha_1, \ldots, \alpha_i\}$, with uncorrelated umbrae similar to an umbra $\alpha$, then $\alpha_B = \alpha^{\cdot |B|}$.

The notation (3.14) can be easily extended to umbral polynomials and dot products as follows:

$$p_M = \prod_{p \in \bar{M}} p^{f(p)}, \qquad (\alpha.p)_M = \prod_{p \in \bar{M}} (\alpha.p)^{f(p)}, \qquad [n.(\chi p)]_M = \prod_{p \in \bar{M}} [n.(\chi p)]^{f(p)},$$

where $p$ are umbral polynomials in $K[A]$ and $\alpha \in A$. For example, using this notation, equivalence (3.9) can be rewritten as

$$[n.(\chi \alpha)]_{P_\lambda} \simeq (n)_{\nu_\lambda} \alpha_\lambda, \tag{3.15}$$

where $P_\lambda = \{\alpha^{(r_1)}, \alpha^{2(r_2)}, \ldots\}$. To avoid misunderstandings, we will specify the multiset $M$ where necessary, in order to know which umbrae occur in $M$.

The notion we are going to introduce is quite natural and allows us to compress and simplify notation.

**Definition 3.1.** A subdivision *of a multiset $M$ is a multiset $S$ of $k \leq |M|$ non-empty submultisets $M_i = (\bar{M}_i, f_i)$ of $M$ such that*

(i) $\bigcup_{i=1}^{k} \bar{M}_i = \bar{M}$;
(ii) $\sum_{i=1}^{k} f_i(\mu) = f(\mu)$ *for every* $\mu \in \bar{M}$.

We note that the notion of subdivision is different from that of multiset partition.

**Example 3.1.** Let $M = \{\mu_1^{(2)}, \mu_2^{(1)}, \mu_3^{(2)}\}$, hence $|M| = 5$. Subdivisions of $M$ are $\{\{\mu_1^{(1)}, \mu_2^{(1)}\}, \{\mu_1^{(1)}\}, \{\mu_3^{(2)}\}\}$ and $\{\{\mu_1^{(1)}\}^{(2)}, \{\mu_2^{(1)}, \mu_3^{(2)}\}\}$.

Let $S = (\bar{S}, g)$ be a subdivision of the multiset $M$. Extending the notation (3.14), we set

$$\mu_S = \prod_{M_i \in \bar{S}} \mu_{M_i}^{g(M_i)} \tag{3.16}$$

and so

$$(n.\mu)_S = \prod_{M_i \in \bar{S}} (n.\mu_{M_i})^{g(M_i)}, \qquad [n.(\chi\mu)]_S = \prod_{M_i \in \bar{S}} [n.(\chi\mu_{M_i})]^{g(M_i)}. \tag{3.17}$$

**Remark 3.1.** A special case of (3.16) is the partition of a set. The notation (3.16) becomes

$$\mu_\pi = \prod_{B \in \pi} \mu_B.$$

We may construct a subdivision of the multiset $M$ in the following steps: assume that the elements of $M$ are all distinct, build a set partition and then replace each element in any block with the original one. Thus, any set partition gives rise to a subdivision. More formally, consider a set $[k]$ of $k$ umbral polynomials. Define the function

$$s : [k] \to \bar{M} \tag{3.18}$$

such that $f(\mu_1)$ elements of $[k]$ go in $\mu_1$, $f(\mu_2)$ elements of $[k]$ go in $\mu_2$ and so on. Now consider a partition $\pi = \{B_1, B_2, \ldots, B_m\}$ of $[k]$ into $m$ blocks. Set

(i) $\bar{M}_i = s(B_i) \subseteq \bar{M}$;
(ii) for any $\mu \in \bar{M}_i$, $f_i(\mu) =$ the number of $p \in B_i$ such that $s(p) = \mu$.

The multiset $S_\pi$, built with $M_i = (\bar{M}_i, f_i)$, $i = 1, 2, \ldots, m$, is the subdivision of $M$ corresponding to the partition $\pi$. Note that $|\pi| = |S_\pi|$. Moreover, it could be $S_{\pi_1} = S_{\pi_2}$ for $\pi_1 \neq \pi_2$.

**Remark 3.2.** If $M$ is a set, it is natural to define $s$ as the identity function so that $S_\pi = \pi$. If $M$ is a subdivision, then $s(p)$ is a multiset.

We simplify equivalences (3.12) and (3.13) by means of the notion of subdivision.

**Proposition 3.8.** *If $S_\pi$ is the subdivision of the multiset $M = \{\alpha^{(k)}\}$ corresponding to a partition $\pi \in \Pi_k$, then*

$$[n.(\chi\alpha)]^k \simeq \sum_{\pi \in \Pi_k} (\chi.\chi)^{\cdot\pi} (n.\alpha)_{S_\pi}, \qquad (n.\alpha)^k \simeq \sum_{\pi \in \Pi_k} [n.(\chi\alpha)]_{S_\pi}. \qquad (3.19)$$

**Proof.** Subdivisions of $M = \{\alpha^{(k)}\}$ are of the type

$$S = \{\underbrace{\{\alpha\}, \ldots, \{\alpha\}}_{r_1}, \underbrace{\{\alpha^{(2)}\}, \ldots, \{\alpha^{(2)}\}}_{r_2}, \ldots\}, \qquad (3.20)$$

with $r_1 + r_2 + \cdots \leq k$ and $r_1 + 2r_2 + \cdots = k$. Via the function $s$ in (3.18), the multiset $S$ corresponds to a partition $\pi$ of $[k]$ with $r_1$ blocks of cardinality 1, $r_2$ blocks of cardinality 2 and so on, so that

$$(n.\alpha)_{S_\pi} = \prod_{M_i \in \bar{S}} (n.\alpha_{M_i})^{g(M_i)} = (n.\alpha)^{r_1} (n.\alpha^2)^{r_2} \cdots,$$

by which the former equivalence in (3.19) is the result of (3.13). The latter equivalence follows from (3.12) by similar arguments. □

A natural extension of the notation to umbral polynomials leads to the following corollary.

**Corollary 3.3.** *If $M$ is a multiset of umbral polynomials, then*

$$[n.(\chi p)]_M \simeq \sum_{\pi \in \Pi_k} (\chi.\chi)^{\cdot\pi} (n.p)_{S_\pi}, \qquad (n.p)_M \simeq \sum_{\pi \in \Pi_k} [n.(\chi p)]_{S_\pi}, \qquad (3.21)$$

*where $k$ is the length of $M$ and $S_\pi$ is the subdivision corresponding to the partition $\pi \in \Pi_k$.*

## 4. *k*-statistics and polykays

If $a_1, a_2, \ldots$ are moments of a random variable and $\kappa_1, \kappa_2, \ldots$ are its cumulants, then

$$\kappa_i = \sum_{\lambda \vdash i} (-1)^{\nu_\lambda - 1} (\nu_\lambda - 1)! \, d_\lambda a_\lambda, \qquad (4.1)$$

where $\lambda = (1^{r_1}, 2^{r_2}, \ldots), a_\lambda = a_1^{r_1} a_2^{r_2} \cdots$ and $d_\lambda$ is given in (2.2). If the umbra $\alpha$ represents the sequence $a_1, a_2, \ldots$, then the sequence of its cumulants $\kappa_1, \kappa_2, \ldots$ is represented by the $\alpha$-cumulant umbra $\chi.\alpha$, as follows by comparing (2.8) and (4.1) (for more details, see Di Nardo and Senato (2006a)). The $i$th $k$-statistic $k_i$ is the unique symmetric unbiased estimator of the cumulant $\kappa_i$ of a given statistical distribution, that is,

$$E[k_i] = \kappa_i$$

(see Stuart and Ord (1987)). $k$-statistics are usually expressed in terms of sums of the $r$th powers of the data points

$$s_r = \sum_{i=1}^{n} X_i^r.$$

In the following, we will give an umbral expression of $k$-statistics by using umbral power sum symmetric polynomials in $n$ uncorrelated and similar umbrae, that is, $n.\alpha^r$.

**Theorem 4.1 ($k$-statistics).** *For $i \leq n$, we have*

$$(\chi.\alpha)^i \simeq \sum_{\lambda \vdash i} \frac{(\chi.\chi)^{\nu_\lambda}}{(n)_{\nu_\lambda}} d_\lambda \sum_{\pi \in \Pi_{\nu_\lambda}} (\chi.\chi)^{\cdot \pi} (n.\alpha)_{S_\pi}, \tag{4.2}$$

*where $\lambda = (1^{r_1}, 2^{r_2}, \ldots)$ runs over all partitions of the integer $i$ and $S_\pi$ is the subdivision of the multiset $P_\lambda = \{\alpha^{(r_1)}, \alpha^{2(r_2)}, \ldots\}$ corresponding to the partition $\pi \in \Pi_{\nu_\lambda}$.*

**Proof.** By replacing equivalence (3.15) in (2.8), we have

$$(\chi.\alpha)^i \simeq \sum_{\lambda \vdash i} \frac{(\chi.\chi)^{\nu_\lambda}}{(n)_{\nu_\lambda}} d_\lambda [n.(\chi\alpha)]_{P_\lambda},$$

where $P_\lambda = \{\alpha^{(r_1)}, \alpha^{2(r_2)}, \ldots\}$. Equivalence (4.2) is the result of the former in (3.21), where $p$ has been replaced by $\alpha$ and the multiset $M$ by $P_\lambda$. □

Since $E[(\chi.\alpha)^i] = \kappa_i$, equivalence (4.2) gives the umbral expression of the $i$th cumulant in terms of umbral power sum symmetric polynomials, that is, the $i$th $k$-statistic.

**Example 4.1.** The partitions of the integer 3 are $\{(1^3), (1^1, 2^1), (3^1)\}$, of length $3, 2, 1$, respectively. Hence,

$$x_{(\nu_\lambda)} = \begin{cases} 2, & \text{for } \lambda = (1^3), \\ -1, & \text{for } \lambda = (1^1, 2^1), \\ 1, & \text{for } \lambda = (3^1), \end{cases} \qquad d_\lambda = \begin{cases} 1, & \text{for } \lambda = (1^3), \\ 3, & \text{for } \lambda = (1^1, 2^1), \\ 1, & \text{for } \lambda = (3^1). \end{cases}$$

From (4.2) and Table 1, we have

$$k_3 = \frac{2}{(n)_3}(2s_3 - 3s_1s_2 + s_1^3) - \frac{3}{(n)_2}(-s_3 + s_1s_2) + \frac{1}{n}s_3$$
$$= \frac{n^2 s_3 - 3n s_1 s_2 + 2s_1^3}{(n)_3},$$

which is the well-known *k*-statistic of order 3.

Products of *k*-statistics are known as *polykays*. Indeed, the symmetric statistic $k_{r,\dots,t}$ such that

$$E[k_{r,\dots,t}] = \kappa_r \cdots \kappa_t$$

(where $\kappa_r, \dots, \kappa_t$ are cumulants) generalizes *k*-statistics and these were originally called *generalized k-statistics* by Dressel (1940). Being a product of cumulants, the umbral expression of a polykay is simply

$$k_{r,\dots,t} \simeq (\chi \boldsymbol{.} \alpha)^r \cdots (\chi' \boldsymbol{.} \alpha')^t, \tag{4.3}$$

with $\chi, \dots, \chi'$ being uncorrelated umbrae likewise $\alpha, \dots, \alpha'$ satisfying $\alpha \equiv \cdots \equiv \alpha'$. The following proposition provides the right-hand product of (4.3) in terms of umbral power sum symmetric polynomials.

**Theorem 4.2 (Polykays).** *If $r + \cdots + t \leq n$, then*

$$k_{r,\dots,t} = \sum_{(\lambda \vdash r,\dots,\eta \vdash t)} \frac{(\chi \boldsymbol{.} \chi)^{\nu_\lambda} \cdots (\chi \boldsymbol{.} \chi)^{\nu_\eta} d_\lambda \cdots d_\eta}{(n)_{\nu_\lambda + \cdots + \nu_\eta}} \sum_{\pi \in \Pi_{\nu_\lambda + \cdots + \nu_\eta}} (\chi \boldsymbol{.} \chi)^{\boldsymbol{.} \pi} (n \boldsymbol{.} \alpha)_{S_\pi}, \tag{4.4}$$

*where $\lambda = (1^{r_1}, 2^{r_2}, \dots)$ runs over all partitions of $r$, $\eta = (1^{t_1}, 2^{t_2}, \dots)$ runs over all partitions of $t$ and $S_\pi$ is the subdivision of the multiset*

$$P_{\lambda + \cdots + \eta} = \{\alpha^{(r_1 + \cdots + t_1)}, \alpha^{2(r_2 + \cdots + t_2)}, \dots\}$$

**Table 1.** Formula (4.2)

| $\lambda$ | $P_\lambda$ | $\pi \in \Pi_{\nu_\lambda}$ | $S_\pi$ | $\sharp S_\pi$ | $(\chi \boldsymbol{.} \chi)^{\boldsymbol{.} \pi}(n \boldsymbol{.} \alpha)_{S_\pi}$ |
|---|---|---|---|---|---|
| $(1^3)$ | $\{\alpha^{(3)}\}$ | $\{\{p_1, p_2, p_3\}\}$ | $\{\{\alpha, \alpha, \alpha\}\}$ | 1 | $2(n \boldsymbol{.} \alpha^3)$ |
| | | $\{\{p_1\}, \{p_2, p_3\}\}$ | $\{\{\alpha\}, \{\alpha, \alpha\}\}$ | 3 | $-(n \boldsymbol{.} \alpha)(n \boldsymbol{.} \alpha^2)$ |
| | | $\{\{p_1\}, \{p_2\}, \{p_3\}\}$ | $\{\{\alpha\}, \{\alpha\}, \{\alpha\}\}$ | 1 | $(n \boldsymbol{.} \alpha)^3$ |
| $(1^1, 2^1)$ | $\{\alpha^{(1)}, \alpha^{2(1)}\}$ | $\{\{\alpha, \alpha^2\}\}$ | $\{\{\alpha, \alpha^2\}\}$ | 1 | $-(n \boldsymbol{.} \alpha^3)$ |
| | | $\{\{\alpha\}, \{\alpha^2\}\}$ | $\{\{\alpha\}, \{\alpha^2\}\}$ | 1 | $(n \boldsymbol{.} \alpha)(n \boldsymbol{.} \alpha^2)$ |
| $(3^1)$ | $\{\alpha^3\}$ | $\{\{\alpha^3\}\}$ | $\{\{\alpha^3\}\}$ | 1 | $(n \boldsymbol{.} \alpha^3)$ |

**Table 2.** Formula (4.4)

| $\lambda + \eta$ | $P_{\lambda+\eta}$ | $S_\pi$ | $\sharp S_\pi$ | $(\chi \cdot \chi)^{\cdot \pi}(n.\alpha)_{S_\pi}$ |
|---|---|---|---|---|
| $(2^2)$ | $\{\alpha^2, \alpha^2\}$ | $\{\{\alpha^2, \alpha^2\}\}$ | 1 | $-(n.\alpha^4)$ |
| | | $\{\{\alpha^2\}, \{\alpha^2\}\}$ | 1 | $(n.\alpha^2)^2$ |
| $(1^2, 2)$ | $\{\alpha, \alpha, \alpha^2\}$ | $\{\{\alpha, \alpha, \alpha^2\}\}$ | 1 | $2(n.\alpha^4)$ |
| | | $\{\{\alpha, \alpha\}, \{\alpha^2\}\}$ | 1 | $-(n.\alpha^2)^2$ |
| | | $\{\{\alpha, \alpha^2\}, \{\alpha\}\}$ | 2 | $-(n.\alpha^3)(n.\alpha)$ |
| | | $\{\{\alpha\}, \{\alpha\}, \{\alpha^2\}\}$ | 1 | $(n.\alpha)^2(n.\alpha^2)$ |
| $(1^4)$ | $\{\alpha, \alpha, \alpha, \alpha\}$ | $\{\{\alpha, \alpha, \alpha, \alpha\}\}$ | 1 | $-6(n.\alpha^4)$ |
| | | $\{\{\alpha\}, \{\alpha, \alpha, \alpha\}\}$ | 4 | $2(n.\alpha)(n.\alpha^3)$ |
| | | $\{\{\alpha, \alpha\}, \{\alpha, \alpha\}\}$ | 3 | $(n.\alpha^2)^2$ |
| | | $\{\{\alpha\}, \{\alpha\}, \{\alpha, \alpha\}\}$ | 6 | $-(n.\alpha)^2(n.\alpha^2)$ |
| | | $\{\{\alpha\}, \{\alpha\}, \{\alpha\}, \{\alpha\}\}$ | 1 | $(n.\alpha)^4$ |

*corresponding to the partition* $\pi \in \Pi_{\nu_\lambda + \cdots + \nu_\eta}$.

**Proof.** Note that

$$\alpha_\lambda \cdots \alpha'_\eta \simeq (\alpha_{j_1})^{\cdot(r_1+\cdots+t_1)}(\alpha_{j_2}^2)^{\cdot(r_2+\cdots+t_2)} \cdots \simeq \alpha_{\lambda+\cdots+\eta},$$

where we have denoted by $\lambda + \cdots + \eta$ the integer partition $(1^{r_1+\cdots+t_1}, 2^{r_2+\cdots+t_2}, \ldots)$. Replacing the right-hand product in (4.3) by the product of uncorrelated (2.8) and using (3.15), we have

$$k_{r,\ldots,t} \simeq \sum_{(\lambda \vdash r, \ldots, \eta \vdash t)} \frac{(\chi \cdot \chi)^{\nu_\lambda} \cdots (\chi \cdot \chi)^{\nu_\eta} d_\lambda \cdots d_\eta}{(n)_{\nu_\lambda + \cdots + \nu_\eta}} [n.(\chi\alpha)]_{P_{\lambda+\cdots+\eta}}.$$

Equivalence (4.4) is the result of the former in (3.21), where $p$ has been replaced by $\alpha$ and the multiset $M$ by $P_{\lambda+\cdots+\eta}$. □

*Example 4.2.* Assume $r = t = 2$. In order to express $k_{2,2}$, we need to consider the pairs of partitions $\{((2),(2)); ((1^2),(2)); ((2),(1^2)); ((1^2),(1^2))\}$. Hence,

$$(n)_{\nu_\lambda + \nu_\eta} = \begin{cases} (n)_2, & \text{for } (\lambda, \eta) = ((2),(2)), \\ (n)_3, & \text{for } (\lambda, \eta) = ((1^2),(2)) \quad \text{and} \quad (\lambda, \eta) = ((2),(1^2)), \\ (n)_4, & \text{for } (\lambda, \eta) = ((1^2),(1^2)) \end{cases}$$

and $x_{\nu_{(2)}} = 1, x_{\nu_{(1^2)}} = -1, d_{(2)} = d_{(1^2)} = 1$. From Table 2, we have

$$k_{2,2} = \frac{-s_4 + s_2^2}{(n)_2} - \frac{2(2s_4 - s_2^2 - 2s_3 s_1 + s_1^2 s_2)}{(n)_3} + \frac{-6s_4 + 8s_3 s_1 + 3s_2^2 - 6s_1^2 s_2 + s_1^4}{(n)_4}.$$

## 5. Multivariate moments and multivariate cumulants

In this section, we define multivariate moments and multivariate cumulants of umbral monomials by means of technicalities introduced in Section 3.2. In the following, $M = \{\mu_1^{(f(\mu_1))}, \mu_2^{(f(\mu_2))}, \ldots, \mu_r^{(f(\mu_r))}\}$ denotes a multiset of length $k$.

**Definition 5.1.** *A multivariate moment is the element of $K$ corresponding to the umbral monomial $\mu_M$ via the evaluation $E$, that is,*

$$E[\mu_M] = m_{t_1 \ldots t_r},$$

*where $t_i = f(\mu_i)$ for $i = 1, 2, \ldots, r$.*

For example, if $M = \{\mu_1^{(1)}, \mu_2^{(2)}, \mu_3^{(1)}\}$, we have $E[\mu_M] = m_{121}$. When the umbral monomials $\mu_i$ are uncorrelated, $m_{t_1 \ldots t_r}$ becomes the product of moments of $\mu_i$.

**Definition 5.2.** *A multivariate cumulant is the element of $K$ corresponding to*

$$E[(\chi \boldsymbol{\cdot} \mu)_M] = \kappa_{t_1 \ldots t_r}, \tag{5.1}$$

*where $t_i = f(\mu_i)$ for $i = 1, 2, \ldots, r$.*

Replacing the umbra $\chi$ in (5.1) by the unity umbra $u$, we get multivariate moments. Setting $\mu_i = \alpha$ for $i = 1, 2, \ldots, r$, we have $M = \{\alpha^{(k)}\}$ and $(\chi \boldsymbol{\cdot} \alpha)_M = (\chi \boldsymbol{\cdot} \alpha)^k$, that is, the ordinary $k$th cumulant.

The notion of the *generalized cumulant*, introduced by McCullagh (1984), is translated into the umbral framework as follows:

$$E[(\chi \boldsymbol{\cdot} \mu)_\pi] = \kappa_\pi,$$

where $\pi$ is a partition of a set of umbral monomials. For example, when

$$\pi = \{\{\mu_1, \mu_2\}, \{\mu_3\}, \{\mu_4, \mu_5\}\},$$

we have

$$E[(\chi \boldsymbol{\cdot} \mu)_\pi] = E[(\chi \boldsymbol{\cdot} \mu_1 \mu_2)(\chi \boldsymbol{\cdot} \mu_3)(\chi \boldsymbol{\cdot} \mu_4 \mu_5)] = \kappa^{12,3,45},$$

where $\kappa^{12,3,45}$ is the McCullagh's notation. In particular, when $\pi = \{\{\mu_1\}, \ldots, \{\mu_r\}\}$, we have $(\chi \boldsymbol{\cdot} \mu)_\pi \equiv (\chi \boldsymbol{\cdot} \mu)_{[r]}$, that is, the joint cumulant of $\mu_1, \ldots, \mu_r$.

The next proposition allows us to express multivariate cumulants in terms of multivariate moments. We only need an extension of (3.11), where the partition $\pi$ is replaced by a subdivision $S$.

Assume $S = \{M_1^{(g(M_1))}, M_2^{(g(M_2))}, \ldots, M_j^{(g(M_j))}\}$ to be a subdivision of the multiset $M$. Let

$$\mu^{\boldsymbol{\cdot} S} \equiv (\mu_{M_1})^{\boldsymbol{\cdot} g(M_1)} \cdots (\mu'_{M_j})^{\boldsymbol{\cdot} g(M_j)}, \tag{5.2}$$

where $\mu_{M_t}$ are uncorrelated umbral monomials. If $M = \{\alpha^{(k)}\}$, then subdivisions are of type (3.20) and $\alpha^{\cdot S} \equiv \alpha_\lambda$, with $\lambda = (1^{r_1}, 2^{r_2}, \ldots) \vdash k$.

**Proposition 5.1.** *If $S_\pi$ is the subdivision of the multiset $M$ corresponding to the partition $\pi \in \Pi_k$, then*

$$(\chi \boldsymbol{.} \mu)_M \simeq \sum_{\pi \in \Pi_k} (\chi \boldsymbol{.} \chi)^{|\pi|} \mu^{\cdot S_\pi}. \tag{5.3}$$

**Proof.** By using set partitions instead of integer partitions, equivalence (2.8) can be written as

$$(\chi \boldsymbol{.} \alpha)^k \simeq \sum_{\pi \in \Pi_k} (\chi \boldsymbol{.} \chi)^{|\pi|} \alpha^{\cdot \pi} \simeq \sum_{\pi \in \Pi_k} (\chi \boldsymbol{.} \chi)^{|\pi|} \alpha^{\cdot S},$$

as has already been done for (3.12) and (3.13). We replace $(\chi \boldsymbol{.} \alpha)^k$ by $(\chi \boldsymbol{.} \mu)_N$, where $N = \{\alpha^{(k)}\}$, and $\alpha^{\cdot S}$ by $\mu^{\cdot S_\pi}$, where $S_\pi$ is the subdivision corresponding to the partition $\pi$. We have

$$(\chi \boldsymbol{.} \mu)_N \simeq \sum_{\pi \in \Pi_k} (\chi \boldsymbol{.} \chi)^{|\pi|} \mu^{\cdot S_\pi}$$

and the result holds for a more general multiset $M$ by using umbral substitutions. □

When $M$ is a set of $k$ different elements, Proposition 5.1 gives the following well-known relations among joint cumulants and multivariate moments of a random vector $(X_1, \ldots, X_k)$:

$$\kappa(X_1, \ldots, X_k) = \sum_{\pi \in \Pi_k} (|\pi| - 1)!(-1)^{|\pi|-1} \prod_{B \in \pi} E\left(\prod_{j \in B} X_j\right).$$

Equivalence (5.3) can be inverted in order to express multivariate moments in terms of multivariate cumulants. Indeed, as proven in Di Nardo and Senato (2006a), if $\kappa$ is the $\mu$-cumulant, then $\mu \equiv \beta \boldsymbol{.} \kappa$, where $\beta$ is the Bell umbra. In (5.3), replace the umbra $\chi$ by the umbra $\beta$ and the umbral monomial $\mu_j$ by a more general polynomial $p_j$. Since the factorial moments of the Bell umbra are all equal to 1 (see Example 2.2), we have

$$E[(\beta \boldsymbol{.} p)_M] = \sum_{\pi \in \Pi_k} p^{\cdot S_\pi}.$$

Now replace the umbral polynomial $p$ by the umbral polynomial $\kappa^\mu \equiv \chi \boldsymbol{.} \mu$. From the third equivalence in (2.7), we have $\beta \boldsymbol{.} \chi \boldsymbol{.} \mu \equiv \mu$. This last equivalence allows to prove the following proposition.

**Proposition 5.2.** *If $S_\pi$ is the subdivision of the multiset $M$ corresponding to the partition $\pi \in \Pi_k$, then*

$$\mu_M \simeq \sum_{\pi \in \Pi_k} (\chi \boldsymbol{\cdot} \mu)^{\boldsymbol{\cdot} S_\pi}.$$

The next propositions and corollary are given by Speed (1983) and McCullagh (1984), using different methods and notation.

**Proposition 5.3.** *Assume the set of umbral monomials $\{\mu_1, \mu_2, \ldots, \mu_i\}$ to be the union of two subsets $\{\mu_{j_1}, \ldots, \mu_{j_t}\}$ and $\{\mu_{k_1}, \ldots, \mu_{k_s}\}$, with $s + t = i$, such that umbral monomials belonging to different subsets are uncorrelated. We have*

$$(\chi \boldsymbol{\cdot} \mu)_{[i]} = (\chi \boldsymbol{\cdot} \mu_1) \cdots (\chi \boldsymbol{\cdot} \mu_i) \simeq 0. \tag{5.4}$$

**Proof.** Let $P = \sum_{l=1}^{t} \mu_{j_l}$ and $Q = \sum_{l=1}^{s} \mu_{k_l}$. The polynomials $P$ and $Q$ are uncorrelated. By virtue of (2.9), we have

$$\chi \boldsymbol{\cdot} (P + Q) \equiv \chi \boldsymbol{\cdot} P \dotplus \chi \boldsymbol{\cdot} Q,$$

that is, products involving powers of $\chi \boldsymbol{\cdot} P$ and $\chi \boldsymbol{\cdot} Q$ vanish, (5.4) being a special case. □

When the umbral monomials $\mu_j$ are interpreted as random variables, equivalence (5.4) states a well-known result: if some of the random variables are uncorrelated with all others, then their joint cumulant is zero.

**Corollary 5.1.** *If $\pi$ is a partition of the set $\{\mu_1, \ldots, \mu_i\}$, then*

$$\mu^{\boldsymbol{\cdot} \pi} \simeq \sum_{\substack{\tau \in \Pi_i \\ \tau \leq \pi}} (\chi \boldsymbol{\cdot} \mu)^{\boldsymbol{\cdot} \tau}.$$

**Proof.** Observe that any partition $\tau$ satisfying $\tau > \pi$, has at least one block, say $B$, that is the union of two or more blocks of the partition $\pi$. By Proposition 5.3, we have $E[(\chi \boldsymbol{\cdot} \mu)_B] = 0$ since, in $B$, there is at least one umbral monomial uncorrelated with each other. □

**Proposition 5.4.** *If $\pi$ is a partition of the set $\{\mu_1, \ldots, \mu_i\}$, then*

$$(\chi \boldsymbol{\cdot} \mu)_\pi \simeq \sum_{\substack{\rho \in \Pi_i \\ \rho \geq \pi}} (\chi \boldsymbol{\cdot} \chi)^{|\rho|} \mu^{\boldsymbol{\cdot} \rho}.$$

**Proof.** Assume that the partition $\pi = \{B_1, B_2, \ldots, B_k\}$ has $k \leq i$ blocks. From (5.3), we have

$$(\chi \boldsymbol{\cdot} p)_{[k]} \simeq \sum_{\tau \in \Pi_k} (\chi \boldsymbol{\cdot} \chi)^{|\tau|} p^{\boldsymbol{\cdot} \tau},$$

where $p_j = \mu_{B_j}$ for $j = 1, 2, \ldots, k$. This result follows by observing that $(\chi \cdot p)_{[k]} \simeq (\chi \cdot \mu)_\pi$ and that each partition $\tau \in \Pi_k$ corresponds to a partition $\rho \in \Pi_i$ such that $\rho \geq \tau$ and $p^{\cdot \tau} \simeq \mu^{\cdot \rho}$. □

Corollary 5.1 and Proposition 5.4 give the following result.

**Corollary 5.2.** *If $\pi$ is a partition of the set $\{\mu_1, \ldots, \mu_i\}$, then*

$$(\chi \cdot \mu)_\pi \simeq \sum_{\substack{\rho \in \Pi_i \\ \rho \geq \pi}} (\chi \cdot \chi)^{|\rho|} \sum_{\substack{\tau \in \Pi_i \\ \tau \leq \rho}} (\chi \cdot \mu)^{\cdot \tau}. \tag{5.5}$$

Following the same arguments used by McCullagh (1984), equivalence (5.5) can be rewritten as

$$(\chi \cdot \mu)_\pi \simeq \sum_{\substack{\rho \in \Pi_i \\ \rho \vee \pi = 1}} (\chi \cdot \mu)^{\cdot \rho},$$

where ∨ means the least upper bound and 1 is the full set.

## 6. Multivariate *k*-statistics and multivariate polykays

In order to introduce an umbral version of multivariate $k$-statistics, we must generalize the result of Theorem 3.2.

**Theorem 6.1.** *Let $S = \{M_1^{(g(M_1))}, M_2^{(g(M_2))}, \ldots, M_j^{(g(M_j))}\}$ be a subdivision of $M$. We have*

$$\mu^{\cdot S} \simeq \frac{1}{(n)_{|S|}} [n \cdot (\chi \mu)]_S, \tag{6.1}$$

*where $\mu^{\cdot S}$ has been defined in (5.2).*

**Proof.** Observe that if $M = \{\alpha^{(k)}\}$, then

$$[n \cdot (\chi \alpha)]_S \equiv [n \cdot (\chi \alpha)]^{r_1} [n \cdot (\chi \alpha^2)]^{r_2} \cdots,$$

where $S$ is a subdivision of $M$ of type (3.20). Equivalence (3.15) can then be rewritten as

$$[n \cdot (\chi \alpha)]_S \simeq (n)_{|S|} \alpha^{\cdot S}$$

and, by using umbral substitutions, the result holds for a more general multiset $M$. □

**Table 3.** Formula (6.2)

| $\pi$ | $\tau$ | $S_\tau$ | $(n.\mu)_{S_\tau}$ |
|---|---|---|---|
| $\{\{p_1, p_2, p_3\}\}$ | $\{B_1\}$ | $\{\{\mu_1, \mu_1, \mu_2\}\}$ | $n.(\mu_1^2 \mu_2)$ |
| $\{\{p_1\}, \{p_2, p_3\}\}$ | $\{B_1, B_2\}$ | $\{\{\mu_1, \mu_1, \mu_2\}\}$ | $n.(\mu_1^2 \mu_2)$ |
|  | $\{\{B_1\}, \{B_2\}\}$ | $\{\{\mu_1\}, \{\mu_1, \mu_2\}\}$ | $n.(\mu_1) n.(\mu_1 \mu_2)$ |
| $\{\{p_2\}, \{p_1, p_3\}\}$ | $\{B_1, B_2\}$ | $\{\{\mu_1, \mu_1, \mu_2\}\}$ | $n.(\mu_1^2 \mu_2)$ |
|  | $\{\{B_1\}, \{B_2\}\}$ | $\{\{\mu_1\}, \{\mu_1, \mu_2\}\}$ | $n.(\mu_1) n.(\mu_1 \mu_2)$ |
| $\{\{p_3\}, \{p_1, p_2\}\}$ | $\{B_1, B_2\}$ | $\{\{\mu_2, \mu_1, \mu_1\}\}$ | $n.(\mu_1^2 \mu_2)$ |
|  | $\{\{B_1\}, \{B_2\}\}$ | $\{\{\mu_2\}, \{\mu_1, \mu_1\}\}$ | $n.(\mu_2) n.(\mu_1^2)$ |
| $\{\{p_1\}, \{p_2\}, \{p_3\}\}$ | $\{B_1, B_2, B_3\}$ | $\{\{\mu_1, \mu_1, \mu_2\}\}$ | $n.(\mu_1^2 \mu_2)$ |
|  | $\{\{B_1\}, \{B_2, B_3\}\}$ | $\{\{\mu_1\}, \{\mu_1, \mu_2\}\}$ | $n.(\mu_1) n.(\mu_1 \mu_2)$ |
|  | $\{\{B_2\}, \{B_1, B_3\}\}$ | $\{\{\mu_1\}, \{\mu_1, \mu_2\}\}$ | $n.(\mu_1) n.(\mu_1 \mu_2)$ |
|  | $\{\{B_3\}, \{B_1, B_2\}\}$ | $\{\{\mu_2\}, \{\mu_1, \mu_1\}\}$ | $n.(\mu_2) n.(\mu_1^2)$ |
|  | $\{\{B_1\}, \{B_2\}, \{B_3\}\}$ | $\{\{\mu_1\}, \{\mu_1\}, \{\mu_2\}\}$ | $n.(\mu_2)[n.(\mu_1)]^2$ |

**Theorem 6.2 (Multivariate *k*-statistics).** *If $n > |M| = k$, then*

$$(\chi.\mu)_M \simeq \sum_{\pi \in \Pi_k} \frac{(\chi.\chi)^{|\pi|}}{(n)_{|\pi|}} \sum_{\tau \in \Pi_{|\pi|}} (\chi.\chi)^{.\tau} (n.\mu)_{S_\tau}, \tag{6.2}$$

*where $S_\tau$ is the subdivision of $M$ corresponding to the partition $\tau$ of the set built with the blocks of $\pi$.*

**Proof.** By Theorem 6.1 and Proposition 5.1, we have

$$(\chi.\mu)_M \simeq \sum_{\pi \in \Pi_k} \frac{(\chi.\chi)^{|\pi|}}{(n)_{|S_\pi|}} [n.(\chi\mu)]_{S_\pi}.$$

The result follows from the first part of (3.21) and by recalling that $|\pi| = |S_\pi|$. □

It is interesting to note the similarity between expressions (4.2) and (6.2). On the right-hand side of (6.2), the set partition replaces the integer partition, the support of $M$ having a cardinality greater than 1.

***Example 6.1.*** In order to express $k_{21}$, take the multiset $M = \{\mu_1^{(2)}, \mu_2^{(1)}\}$ of length $k = 3$. Define the function $s : [3] \to \{\mu_1, \mu_2\}$ such that $s(p_1) = s(p_2) = \mu_1$ and $s(p_3) = \mu_2$. Let $s_{p,q} \simeq n.(\mu_1^p \mu_2^q)$. From Table 3, we have

$$k_{21} \simeq (\chi.\mu_1)^2 (\chi.\mu_2) \simeq \frac{1}{(n)_3} [n^2 s_{2,1} - 2n s_{1,0} s_{1,1} - n s_{2,0} s_{0,1} + 2 s_{1,0}^2 s_{0,1}].$$

Multivariate polykays were introduced by Robson (1957). The symmetric statistic $k_{t_1...t_r,...,l_1...l_m}$ such that

$$E[k_{t_1...t_r,...,l_1...l_m}] = \kappa_{t_1...t_r} \cdots \kappa_{l_1...l_m}$$

(where $\kappa_{t_1...t_r}$ is a multivariate cumulant) generalizes polykays. Being a product of uncorrelated multivariate cumulants, the umbral expression of a multivariate polykay is simply

$$k_{t_1...t_r,...,l_1...l_m} \simeq (\chi \boldsymbol{\cdot} \mu)_T \cdots (\chi' \boldsymbol{\cdot} \mu')_L,$$

where $\chi$ and $\chi'$ are uncorrelated likewise the umbral monomials $\mu \in T$ and $\mu' \in L$ where

$$T = \{\mu_1^{(t_1)}, \ldots, \mu_r^{(t_r)}\}, \ldots, L = \{\mu'^{(l_1)}_1, \ldots, \mu'^{(l_m)}_m\}.$$

Let $S$ be a subdivision of $T$ and $S^*$ a subdivision of $L$.

**Theorem 6.3 (Multivariate polykays).** *For $n > |T| + \cdots + |L|$, we have*

$$k_{t_1...t_r,...,l_1...l_m} \simeq \sum_{(\pi \in \Pi_{|T|}, \ldots, \tilde{\pi} \in \Pi_{|L|})} \frac{(\chi \boldsymbol{\cdot} \chi)^{|\pi|} \cdots (\chi' \boldsymbol{\cdot} \chi')^{|\tilde{\pi}|}}{(n)_{|\pi|+\cdots+|\tilde{\pi}|}} \sum_{\tau \in \Pi_{|\pi|+\cdots+|\tilde{\pi}|}} (\chi \boldsymbol{\cdot} \chi)^{\boldsymbol{\cdot} \tau} (n \boldsymbol{\cdot} p)_{S_\tau}, \tag{6.3}$$

*where $S_\tau$ is the subdivision of the multiset obtained by the disjoint union of $T, \ldots, L$, with no uncorrelation labels and corresponding to the partition $\tau$ of the set built with the blocks of $\{\pi, \ldots, \tilde{\pi}\}$.*

**Proof.** Observe that if $S_\pi$ is the subdivision of $T$ corresponding to the partition $\pi \in \Pi_{|T|}$ and $S_{\tilde{\pi}}$ is the subdivision of $L$ corresponding to the partition $\tilde{\pi} \in \Pi_{|L|}$, then we have

$$\mu^{\boldsymbol{\cdot} S_\pi} \cdots \nu^{\boldsymbol{\cdot} S_{\tilde{\pi}}} \simeq \mu_{S_\pi + \cdots + S_{\tilde{\pi}}},$$

where we have denoted by $S_\pi + \cdots + S_{\tilde{\pi}}$ the subdivision obtained by placing side by side the blocks of subdivisions with no uncorrelated elements. By equivalences (5.3) and (6.1), we have

$$k_{t_1...t_r,...,l_1...l_m} \simeq \sum_{(\pi \in \Pi_{|T|}, \ldots, \tilde{\pi} \in \Pi_{|L|})} \frac{(\chi \boldsymbol{\cdot} \chi)^{|\pi|} \cdots (\chi' \boldsymbol{\cdot} \chi')^{|\tilde{\pi}|}}{(n)_{|\pi|+\cdots+|\tilde{\pi}|}} [n \boldsymbol{\cdot} (\chi p)]_{S_\pi + \cdots + S_{\tilde{\pi}}},$$

where $p$ is the generic element of $S_\pi + \cdots + S_{\tilde{\pi}}$. The result follows by applying the first part of (3.21). □

Also in this case, note the similarity between expressions (4.4) and (6.3). The sum of subdivisions in the second case corresponds to the sum of integer partitions.

**Table 4.** Formula (6.3)

| $\{\pi, \tilde{\pi}\}$ | $\tau$ | $S_\tau$ | $(n.p)_{S'_\tau}$ |
|---|---|---|---|
| $\{\{\mu_1, \mu_2\}, \{\mu'_1\}\}$ | $\{B_1, B_2\}$ | $\{\{\mu_1, \mu_1, \mu_2\}\}$ | $n.(\mu_1^2 \mu_2)$ |
|  | $\{\{B_1\}, \{B_2\}\}$ | $\{\{\mu_1, \mu_2\}, \{\mu_1\}\}$ | $n.(\mu_1 \mu_2) n.\mu_1$ |
| $\{\{\mu_1\}, \{\mu_2\}, \{\mu'_1\}\}$ | $\{\{B_1, B_2, B_3\}\}$ | $\{\{\mu_1, \mu_1, \mu_2\}\}$ | $n.(\mu_1^2 \mu_2)$ |
|  | $\{\{B_1, B_2\}, \{B_3\}\}$ | $\{\{\mu_1, \mu_1\}, \{\mu_2\}\}$ | $n.(\mu_1^2) n.\mu_2$ |
|  | $\{\{B_1, B_3\}, \{B_2\}\}$ | $\{\{\mu_1, \mu_2\}, \{\mu_1\}\}$ | $n.(\mu_1 \mu_2) n.\mu_1$ |
|  | $\{\{B_2, B_3\}, \{B_1\}\}$ | $\{\{\mu_1, \mu_2\}, \{\mu_1\}\}$ | $n.(\mu_1 \mu_2) n.\mu_1$ |
|  | $\{\{B_1\}, \{B_2\}, \{B_3\}\}$ | $\{\{\mu_1\}, \{\mu_2\}, \{\mu_1\}\}$ | $(n.\mu_1)^2 n.\mu_2$ |

**Example 6.2.** In order to express $k_{11,1}$, let $T = \{\mu_1, \mu_2\}, L = \{\mu'_1\}$ and $s_{p,q} \simeq n.(\mu_1^p \mu_2^q)$. From Table 4, we have

$$k_{11,1} \simeq \frac{1}{(n)_2}[s_{1,0} s_{1,1} - s_{2,1}] - \frac{1}{(n)_3}[s_{1,0}^2 s_{0,1} - 2 s_{1,0} s_{1,1} + 2 s_{2,1} - s_{2,0} s_{0,1}].$$

## 7. A fast algorithm for *k*-statistics

It is possible to build a very fast algorithm for *k*-statistics by forfeiting the elegant idea of producing only one algorithm for the whole subject. We are going to prove that *k*-statistics can be recovered through cumulants of compound Poisson random variables. This connection allows us to insert exponential polynomials in formula (4.2), eliminating set partitions.

In order to construct such an algorithm, we need ratios of umbrae. Therefore, we introduce the notion of the *multiplicative inverse* of an umbra. Two umbrae are said to be *multiplicatively inverse* to each other when

$$\alpha \gamma \equiv u. \tag{7.1}$$

In dealing with a saturated umbral calculus, the multiplicative inverse of an umbra is not unique, but any two multiplicatively inverse umbrae of the umbra $\alpha$ are similar. From (7.1), we have

$$a_n g_n = 1 \quad \text{for every } n = 0, 1, 2, \ldots, \text{ that is, } g_n = \frac{1}{a_n},$$

where $a_n$ and $g_n$ are moments of $\alpha$ and $\gamma$, respectively. In the following, a multiplicative inverse of an umbra $\alpha$ will be denoted by $1/\alpha$.

### 7.1. Cumulants of compound Poisson random variables

The main result of this section is that cumulants of an umbra $\alpha$ can be expressed via cumulants of compound Poisson random variables.

Let us consider the umbra $\chi.y.\beta$ introduced in Di Nardo and Senato (2006a), where $y$ is an indeterminate. As well as the umbra $\chi.\beta$ has moments all equal to 1, the moments of this umbra are all equal to $y$, as the following lemma states.

**Lemma 7.1.** *If $\chi$ is the singleton umbra and $\beta$ is the Bell umbra, then*

$$(\chi.y.\beta)^i \simeq y, \qquad i = 1, 2, \ldots. \tag{7.2}$$

**Proof.** Observe that $\chi.y.\beta \equiv \chi.y.\beta.u \equiv (\chi.y).\beta.u$, where the last equivalence follows from the associative law. So, by virtue of (2.5), we have

$$[(\chi.y).\beta]^i \simeq [(\chi.y).\beta.u]^i \simeq \sum_{\lambda \vdash i} (\chi.y)^{\nu_\lambda} d_\lambda u_\lambda. \tag{7.3}$$

On the other hand, $(\chi.y) \equiv \chi y$ and the only power of $(\chi.y)$ different from zero is the one corresponding to $\nu_\lambda = 1$ for which $\lambda = (i)$. Therefore, the sum in (7.3) reduces to $y$. $\square$

As stated in Example 2.4, the umbra $(\chi.y.\beta).\alpha \equiv \chi.(y.\beta.\alpha)$ is the cumulant umbra of a polynomial $\alpha$-partition umbra and corresponds to a compound Poisson random variable of parameter $y$ (Di Nardo and Senato (2001)). Such an umbra is the keystone for building cumulants.

First, let us compute the moments of $n.\chi.(y.\beta.\alpha)$, that is, the sum of $n$ uncorrelated cumulant umbrae of a polynomial $\alpha$-partition umbra.

**Proposition 7.1.** *The polynomial umbra $n.\chi.(y.\beta.\alpha)$ has moments*

$$c_i(y) = \sum_{\lambda \vdash i} y^{\nu_\lambda} (n)_{\nu_\lambda} d_\lambda a_\lambda, \tag{7.4}$$

*with $c_0(y) = 1$ and $c_i(y)$ of degree $i$ for every integer $i$.*

**Proof.** Due to the associative law, we have $n.\chi.(y.\beta.\alpha) \equiv (n.\chi.y).\beta.\alpha$ so that

$$[(n.\chi.y).\beta.\alpha]^i \simeq \sum_{\lambda \vdash i} (n.\chi.y)^{\nu_\lambda} d_\lambda \alpha_\lambda \tag{7.5}$$

is the result of (2.5). Let $c_i(y) = E[(n.\chi.y.\beta.\alpha)^i]$. From (2.1), we have

$$(n.\chi.y)^j \simeq [n.(\chi y)]^j \simeq \sum_{\lambda \vdash j} (n)_{\nu_\lambda} d_\lambda (\chi y)_\lambda.$$

Since
$$(\chi y)_\lambda \simeq \begin{cases} y^j, & \text{for } \lambda = (1^j), \\ 0, & \text{otherwise}, \end{cases}$$
we have $(n.\chi.y)^{\nu_\lambda} \simeq (n)_{\nu_\lambda} y^{\nu_\lambda}$. Replacing this last equivalence in (7.5), the result follows. □

The next theorem is the key to producing a fast algorithm for *k*-statistics. It states that the $\alpha$-cumulant umbra has moments umbrally equivalent to umbral polynomials obtained by replacing $y$ with the umbra $\chi.\chi/n.\chi$ in $c_i(y), i = 1, 2, \ldots$.

**Theorem 7.1.** *If $c_i(y) = E[(n.\chi.y.\beta.\alpha)^i]$, then*
$$(\chi.\alpha)^i \simeq c_i\left(\frac{\chi.\chi}{n.\chi}\right), \qquad i = 1, 2, \ldots. \tag{7.6}$$

**Proof.** In (7.4), replace $y$ by the umbra $\chi.\chi/n.\chi$. We have
$$c_i\left(\frac{\chi.\chi}{n.\chi}\right) \simeq \sum_{\lambda \vdash i} \left(\frac{\chi.\chi}{n.\chi}\right)^{\nu_\lambda} (n)_{\nu_\lambda} d_\lambda \alpha_\lambda. \tag{7.7}$$

Due to the uncorrelation property, we have
$$E\left[c_i\left(\frac{\chi.\chi}{n.\chi}\right)\right] = \sum_{\lambda \vdash i} \frac{E[(\chi.\chi)^{\nu_\lambda}]}{E[(n.\chi)^{\nu_\lambda}]} (n)_{\nu_\lambda} d_\lambda a_\lambda.$$

Recalling that $E[(\chi.\chi)^{\nu_\lambda}] = (-1)^{\nu_\lambda}(\nu_\lambda - 1)!$ and $E[(n.\chi)^{\nu_\lambda}] = (n)_{\nu_\lambda}$, the result follows. □

### 7.2. *k*-statistics via compound Poisson random variables

In the previous section, we have stated the connection between cumulants of an umbra and those of compound Poisson random variables. In order to recover the umbral expressions of *k*-statistics given in (4.2), it is sufficient to express the moments of $n.\chi.x.\beta.\alpha$ in terms of power sums $n.\alpha^i$. This task has been partially accomplished in Section 3.1 for the umbra $n.(\chi\alpha)$. We are going to extend such relations to a more general umbra $n.(\gamma\alpha)$.

**Lemma 7.2.** *In $K[x_1, x_2, \ldots, x_n][A]$, we have*
$$\chi.(\gamma_1 x_1 + \cdots + \gamma_n x_n) \equiv (\chi.\gamma)\sigma, \tag{7.8}$$
*where $\sigma$ is the power sum polynomial umbra and $\{\gamma_i\}_{i=1}^n$ are uncorrelated umbrae similar to an umbra $\gamma$.*

**Proof.** Due to property (2.9), we have

$$\chi.(\gamma_1 x_1 + \cdots + \gamma_n x_n) \equiv \dot{+}_{i=1}^{n}[\chi.(\gamma_i x_i)] \equiv \dot{+}_{i=1}^{n}(\chi.\gamma_i) x_i,$$

where the last equivalence follows due to the fact that $\chi.(c\gamma) \equiv c(\chi.\gamma)$ for every $c \in K$. Since the umbrae $\gamma_i$ are similar to the umbra $\gamma$, we have

$$\dot{+}_{i=1}^{n}(\chi.\gamma_i) x_i \equiv (\chi.\gamma)[\dot{+}_{i=1}^{n} u x_i],$$

from which (7.8) follows. □

**Corollary 7.1.** *If $\sigma$ is the power sum polynomial umbra and $\{\gamma_i\}_{i=1}^{n}$ are uncorrelated umbrae similar to the umbra $\gamma$, then*

$$(\gamma_1 x_1 + \cdots + \gamma_n x_n)^i \simeq \sum_{\lambda \vdash i} d_\lambda (\chi.\gamma)_\lambda \sigma_\lambda.$$

**Proof.** Taking the left dot product with $\beta$ in (7.8) and recalling that $\beta.\chi \equiv u$, we have

$$(\gamma_1 x_1 + \cdots + \gamma_n x_n) \equiv \beta.[(\chi.\gamma)\sigma]. \tag{7.9}$$

The result follows via (2.4). □

Replacing the indeterminate $x_i$ by $\alpha_i$, the next theorem follows immediately.

**Theorem 7.2.** *If $\alpha, \gamma \in A$, then*

$$[n.(\gamma\alpha)]^i \simeq \sum_{\lambda \vdash i} d_\lambda (\chi.\gamma)_\lambda (n.\alpha)^{r_1} (n.\alpha^2)^{r_2} \cdots,$$

*where $\lambda = (1^{r_1}, 2^{r_2}, \ldots) \vdash i$.*

Since $\chi.y.\beta.\alpha \equiv (\chi.y.\beta)\alpha$ (Di Nardo and Senato (2006a), formula (31)), Theorem 7.2 allows us to express the polynomials $c_i(y)$ in terms of power sums $n.\alpha^i$. This is the starting point to prove the following result.

**Theorem 7.3.** *Let*

$$p_n(y) = \sum_{k=1}^{n} y^k S(n,k)(-1)^{k-1}(k-1)!, \tag{7.10}$$

*where $S(n,k)$ are the Stirling numbers of second type. We have*

$$(\chi.\alpha)^i \simeq \sum_{\lambda \vdash i} d_\lambda p_\lambda \left(\frac{\chi.\chi}{n.\chi}\right) (n.\alpha)^{r_1} (n.\alpha^2)^{r_2} \cdots,$$

*where $\lambda = (1^{r_1}, 2^{r_2}, \ldots) \vdash i$ and $p_\lambda(y) = [p_1(y)]^{r_1} [p_2(y)]^{r_2} \cdots$.*

**Proof.** From Theorem 7.2, we have

$$c_i(y) \simeq \sum_{\lambda \vdash i} d_\lambda (\chi \boldsymbol{.} \chi \boldsymbol{.} y \boldsymbol{.} \beta)_\lambda (n \boldsymbol{.} \alpha)^{r_1} (n \boldsymbol{.} \alpha^2)^{r_2} \cdots, \tag{7.11}$$

with $c_i(y)$ the $i$th moment of $n\boldsymbol{.}[(\chi\boldsymbol{.}y\boldsymbol{.}\beta)\alpha]$. Note that $\chi\boldsymbol{.}\chi\boldsymbol{.}y\boldsymbol{.}\beta \equiv u^{<-1>}y\boldsymbol{.}\beta$, so powers of $\chi\boldsymbol{.}\chi\boldsymbol{.}y$ are umbrally equivalent to the exponential umbral polynomials

$$\phi_n(\gamma) = \sum_{k=0}^{n} S(n,k)\gamma^k$$

(Di Nardo and Senato (2006a)) with $\gamma$ replaced by $u^{<-1>}y$, that is,

$$\phi_n(u^{<-1>}y) = \sum_{k=1}^{n} y^k S(n,k)(u^{<-1>})^k.$$

Observing that $p_n(y) = E[\phi_n(u^{<-1>}y)]$, the result follows from (7.11) since $(\chi\boldsymbol{.}\chi\boldsymbol{.}y\boldsymbol{.}\beta)_\lambda \simeq p_\lambda(y)$. □

The MAPLE algorithm, which implements the result of Theorem 7.3, is the following:
```
makeTab := proc(N)
  [seq([N!/mul((x!)^numboccur(y,x) * numboccur(y,x)!,
      x = {op(y)}) * mul(k[i], i = y), mul(S[i], i = y)], y = combinat
      ['partition'](N))]; end :
makeK := proc(N)
  [seq(k[i] = add(combinat['stirling2'](i,j) * x^j * (-1)^(j-1) * (j-1)!,
      j = 1..i), i = 1..N)]; end :
fd := proc(j, h)
  expand(mul(n - t - h, t = 0..j - h)); end :
kstat := proc(N)
  local u, v;
  v := expand(eval(makeTab(N), makeK(N)));
  u := [seq(x^i = (-1)^(i-1) * (i-1)! * fd(N-1, i), i = 1..N)];
  v := expand(eval(v, u));
          1/mul((n - x), x = 0..N - 1) * add(x[1] * x[2], x = v); end :
```

## 8. Concluding remarks

Umbral formulae for *k*-statistics and polykays, either in single or multivariate cases, share a common algorithm to construct multiset subdivisions. When the multiset has the form $\{\alpha^{(i)}\}$, this corresponds to computing partitions of the integer $i$, but integer partitions cannot be employed when the multiset has a more complex form. Even though the idea of constructing the function $s$ is fundamental in constructing umbral formulae involving multisets, this is not efficient from a computational point of view. Indeed, examples have shown how subdivisions may occur more than once in the same formula. We have therefore constructed an algorithm that generates only different subdivisions and enumerates how many times each subdivision occurs; see Di Nardo *et al.* (2008). Such an algorithm is the heart of the procedure Polykays which produces $k$-statistics, polykays, multivariate $k$-statistics and multivariate polykays using less computational time than those implemented by Andrews and Stafford (2000). For single and multivariate $k$-statistics Polykays has computational times comparable with MATHSTATICA (Rose and Smith (2002)) up to order 7 and a little worse for higher orders. Moreover, Polykays allows the computation of multivariate polykays, unlike MATHSTATICA.

**Table 5.** Comparisons between computational times for single $k$-statistics obtained by using the MATHEMATICA procedures of Andrews and Stafford, those of MATHSTATICA and the algorithm constructed via Theorem 7.3 (runs on PC 2.08 GHz, 512MB RAM)

| $i$ | Andrews and Stafford | MATHSTATICA | Fast umbral algorithm |
| --- | --- | --- | --- |
| 8  | 0.72     | 0.03    | 0.00 |
| 10 | 3.49     | 0.08    | 0.00 |
| 12 | 24.80    | 0.20    | 0.00 |
| 14 | 396.34   | 0.56    | 0.05 |
| 16 | 58002.60 | 1.69    | 0.11 |
| 18 | –        | 5.42    | 0.20 |
| 20 | –        | 19.11   | 0.41 |
| 22 | –        | 69.66   | 0.81 |
| 24 | –        | 285.58  | 1.66 |
| 26 | –        | 1551.48 | 3.49 |
| 28 | –        | 6324.28 | 7.78 |

Finally, in Section 7 of this paper, we have proven that the umbral techniques not only provide a unifying structure for the whole subject, but also a new way of improving the computational generation of such estimators. An immediate example is the formula stated in Theorem 7.3, whose derived algorithm realizes the amazing computational times shown in Table 5.

# References

Andrews, D.F. (2001). Asymptotic expansions of moments and cumulants. *Stat. Comput.* **11** 7–16. MR1837140

Andrews, D.F. and Stafford, J.E. (1998). Iterated full partitions. *Stat. Comput.* **8** 189–192.

Andrews, D.F. and Stafford, J.E. (2000). *Symbolic Computation for Statistical Inference. Oxford Statistical Science Series* **21**. Oxford Univ. Press. MR1857192

Di Bucchianico, A. (1997). *Probabilistic and Analytical Aspects of the Umbral Calculus. CWI Tract* **119**. Stichting Mathematisch Centrum Centrum voor Wiskunde en Informatica, Amsterdam. MR1431509

Di Nardo, E. and Senato, D. (2001). Umbral nature of the Poisson random variables. In *Algebraic Combinatorics and Computer Science* (H. Crapo and D. Senato, eds.) 245–266. Milano: Springer. MR1854481

Di Nardo, E. and Senato, D. (2006a). An umbral setting for cumulants and factorial moments. *European J. Combin.* **27** 394–413. MR2206475

Di Nardo, E. and Senato, S. (2006b). A symbolic method for $k$-statistics. *Appl. Math. Lett.* **19** 968–975. MR2240494

Di Nardo, E., Guarino, G. and Senato, D. (2008). Maple algorithms for polykays and multivariate polykays. *Adv. in Appl. Statist.* **8** 19–36.

Doubilet, P. (1972). On the foundations of combinatorial theory VII: Symmetric functions through the theory of distribution and occupancy. *Stud. Appl. Math.* **51** 377–396. MR0429577

Dressel, P.L. (1940). Statistical seminvariants and their estimates with particular emphasis on their relation to algebraic invariants. *Ann. Math. Statist.* **11** 33–57. MR0001503

Fisher, R.A. (1929). Moments and product moments of sampling distributions. *Proc. London Math. Soc. (2)* **30** 199–238.

Good, I.J. (1975). A new formula for cumulants. *Math. Proc. Camb. Philos. Soc.* **78** 333–337. MR0386098

Good, I.J. (1977). A new formula for $k$-statistics. *Ann. Statist.* **5** 224–228. MR0423636

McCullagh, P. (1984). Tensor notation and cumulants of polynomials. *Biometrika* **71** 461–476. MR0775392

McCullagh, P. (1987). *Tensor Methods in Statistics*. London: Chapman and Hall. MR0907286

Kaplan, E.L. (1952). Tensor notation and the sampling cumulants of $k$-statistics. *Biometrika* **39** 319–323. MR0051469

Robson, D.S. (1957). Applications of multivariate polykays to the theory of unbiased ratio-type estimation. *J. Amer. Statist. Assoc.* **52** 511–522. MR0092323

Rose, C. and Smith, M.D. (2002). *Mathematical Statistics with Mathematica*. Berlin: Springer. MR1890491

Roman, S.M. and Rota, G.-C. (1978). The umbral calculus. *Adv. in Math.* **27** 95–188. MR0485417

Rota, G.-C. and Taylor, B.D. (1994). The classical umbral calculus. *SIAM J. Math. Anal.* **25** 694–711. MR1266584

Saliani, S. and Senato, D. (2006). Compactly supported wavelets through the classical umbral calculus. *J. Fourier Anal. Appl.* **12** 27–36. MR2215675

Shen, J. (1999). Combinatorics for wavelets: The umbral refinement equation. *Stud. Appl. Math.* **103** 121–147. MR1704938

Speed, T.P. (1983). Cumulants and partition lattices. *Austral. J. Statist.* **25** 378–388. MR0725217

Speed, T.P. (1986a). Cumulants and partition lattices. II. Generalised $k$-statistics. *J. Aust. Math. Soc. Ser. A* **40** 34–53. MR0809724

Speed, T.P. (1986b). Cumulants and partition lattices. III. Multiply-indexed arrays. *J. Aust. Math. Soc. Ser. A* **40** 161–182. MR0817836

Speed, T.P. (1986c). Cumulants and partition lattices. IV. A.s. convergence of generalised $k$-statistics. *J. Aust. Math. Soc. Ser. A* **41** 79–94. MR0846773

Speed, T.P. and Silcock, H.L. (1988a). Cumulants and partition lattices. V. Calculating generalized $k$-statistics. *J. Aust. Math. Soc. Ser. A* **44** 171–196. MR0922603

Speed, T.P. and Silcock, H.L. (1988b). Cumulants and partition lattices. VI. Variances and covariances of mean squares. *J. Aust. Math. Soc. Ser. A* **44** 362–388. MR0929528

Stuart, A. and Ord, J.K. (1987). *Kendall's Advanced Theory of Statistics*. **1**: *Distribution Theory*, 5th ed. London: Charles Griffin and Company Limited.

Thiele, T.N. (1897). *Elementaer Iagttagelseslaere*. Copenhagen: Gyldendal. Reprinted in English (1931). The theory of observations. *Ann. Math. Statist.* **2** 165–308.

Tukey, J.W. (1950). Some sampling simplified. *J. Amer. Statist. Assoc.* **45** 501–519. MR0040624

Tukey, J.W. (1956). Keeping moment-like sampling computations simple. *Ann. Math. Statist.* **27** 37–54. MR0076236

Wishart, J. (1952). Moment coefficients of the $k$-statistics in samples from a finite population. *Biometrika* **39** 1–13. MR0050223